\documentclass{tran-l}

\usepackage{amssymb}

\usepackage{graphicx}

\copyrightinfo{2006}{American Mathematical Society}

\def\proof{\noindent{\bf{Proof.} }}
\def\sqr#1#2{{\vcenter{\hrule height.#2pt
        \hbox{\vrule width.#2pt height#1pt \kern#1pt
                \vrule width.#2pt}
        \hrule height.#2pt}}}
\def\square{\mathchoice\sqr64\sqr64\sqr{4}3\sqr{3}3}
\def\QED{\hfill$\square$}

\def\tratto{\mbox{\rule{2mm}{.2mm}$\;\!$}}

\def\noi{\noindent }

\def\p{{\mathfrak p}}
\def\P{{\mathfrak P}}
\def\Q{{\mathfrak Q}}

\def\n{{\mathfrak n}}
\def\m{{\mathfrak m}}

\def\C{\mathcal C}

\newcommand{\s}{\bigskip}
\newcommand{\ms}{\medskip}
\newcommand{\sss}{\smallskip}
\newcommand{\f}[1]{\ensuremath{\mathfrak{#1}}}
\newcommand{\ol}[1]{\ensuremath{\overline{#1}}}
\newcommand{\ul}[1]{\ensuremath{\underline{#1}}}

\newcommand{\height}{{\rm{ht}}}

\newcommand{\core}[1]{\ensuremath{{\rm{core}}(#1)}}

\newtheorem{theorem}{Theorem}[section]

\newtheorem{corollary}[theorem]{Corollary}
\newtheorem{lemma}[theorem]{Lemma}
\newtheorem{proposition}[theorem]{Proposition}
\newtheorem{Notation and Discussion}[theorem]{Notation and Discussion}
\newtheorem{Assumptions and Discussion}[theorem]{Assumptions and Discussion}
\newtheorem{Assumptions}[theorem]{Assumptions}

 \theoremstyle{definition}

\newtheorem{definition}[theorem]{Definition}

\newtheorem{remark}[theorem]{Remark}

\newtheorem{example}[theorem]{Example}

\numberwithin{equation}{section}

\begin{document}

\title[Annihilators of Components of the Canonical Module]
{\bf Annihilators of Graded Components of the Canonical Module, \\
 and the Core of Standard Graded Algebras}

%    Only \author and \address are required; other information is %    optional.  Remove any unused author tags.

%    author one information
\author[L. Fouli]{Louiza Fouli}
\address{Department of Mathematics, University of Texas, Austin, TX 78712, USA}
\email{lfouli@math.utexas.edu}
\author[C. Polini]{Claudia Polini}
\address{Department of Mathematics, University of Notre Dame, Notre Dame, IN 46556, USA}
\email{cpolini@nd.edu}
\author[B. Ulrich]{Bernd Ulrich}
\address{Department of Mathematics, Purdue University, West Lafayette, IN 47907, USA}
\email{ulrich@math.purdue.edu}

\thanks{The second and third author gratefully acknowledge partial
support from the NSF. The second author was also supported in part
by the NSA. The first and second author thank the Department of
Mathematics of Purdue University for its hospitality}

%   \subjclass is required.
\subjclass[2000]{Primary 13B21; Secondary 13A30, 13B22, 13C40}

%\date{}

%\dedicatory{}

%    Abstract is required.
\begin{abstract}
We relate the annihilators of graded components of the canonical
module of a graded Cohen-Macaulay ring to colon ideals of powers of the
homogeneous maximal ideal. In particular, we connect
them to the core of the maximal ideal. An application of
our results characterizes Cayley-Bacharach sets of points in terms
of the structure of the core of the maximal ideal of their
homogeneous coordinate ring. In particular, we show that a scheme is
Cayley-Bacharach if and only if the core is a power of the maximal
ideal.
\end{abstract}

\maketitle


\section{Introduction}

\sss

This paper started as an attempt to understand the core of powers of
the homogeneous maximal ideal of a standard graded algebra. The {\it
core} of an arbitrary ideal $I$ in a Noetherian ring is defined as
the intersection of all reductions of $I$, equivalently, of all
ideals over which $I$ is integral. The definition simplifies when
$I$ is a power $\m^n$ of the homogeneous maximal ideal $\m$ of a
standard graded Cohen-Macaulay algebra $R$ of positive dimension
over an infinite field $k$. In this case $\core{\m^n}$ is the
intersection of all ideals generated by systems of parameters
consisting of forms of degree $n$; in fact, it suffices to intersect
finitely many parameter ideals generated by {\it general} forms of
degree $n$ \cite{CPU01, PUV}. One has explicit formulas in this case
that express the core of $\m^n$ as a colon ideal and that are valid
in any characteristic. If ${\rm char} \, k =0$ or if $R$ is
geometrically reduced, then
\[\core{\m^n} = J^{j+(n-1)d+1}:\m^{j} \/ ,\]
where $J$ is an ideal generated by a linear system of parameters of
$R$, $\/j\,$ is any integer $\geq a + d$, and $a$, $d$ denote the
$a$-invariant and the dimension of $R$, respectively (see
Theorem~\ref{formula} and, for prior results, \cite{HySm2, FPU}).
Thus the question arises of what more can be said about the colon
ideals $J^{j+(n-1)d+1}:\m^{j}$ and, most notably, when the obvious
inclusion $ \m^{a+nd+1} \subset J^{j+(n-1)d+1}:\m^{j}$ is an
equality. We were also wondering in what sense the shape of the core
of $\m$ or its powers reflects the geometry of ${\rm Proj}(R)$ as a
subscheme of projective space. As it turns out, the crucial device
in approaching these questions is the graded canonical module
$\omega$ of $R$ and the faithfulness of the submodules generated by
its graded components. Thus the main goal of this paper is to study,
quite generally, the interplay between annihilators of graded
components of $\omega$ on the one hand and colon ideals of powers of
$\m$ on the other hand.

More generally, let $R$ be a standard graded Cohen-Macaulay algebra
over a field $k$, of dimension $d >0$, and write $\m$ for its
homogeneous maximal ideal and $\omega = \omega_R$ for its graded
canonical module. By $a=a(R)$ we denote the $a$-{\it invariant} of
$R$, which is the negative of the initial degree of $\omega$. Recall
that if $J$ is any ideal generated by a linear system of parameters
then `$J$ is a reduction of $\m$ with reduction number
$r_J(\m)=a+d$', which simply means that $\m^i=J^{i-j}\m^j$ for every
$i \geq j \geq a+d$, but for no smaller $j$. It is easy to see that
$J^i:\m^j={\m}^{i-j+a+d}$ for every $i$ and $j \geq
a+d$, provided $R$ is Gorenstein or, more generally, {\it level},
which means that $\omega$ is generated by homogeneous elements of
the same degree, $\omega = [\omega]_{-a}R$. Under these
assumptions, the $R$-submodules $[\omega]_tR$ of $\omega$ generated
by the homogeneous elements of a fixed degree $t$ are all faithful as
long as $t \geq -a$; this is obvious since $\omega=[\omega]_{-a}R$
is faithful and there exists a form of positive degree regular on
$\omega$. The same holds if $R$ is a domain because a suitable shift
of $\omega$ embeds into $R$. Without the additional assumptions on
the ring neither the statement about the colon ideals nor the one
about the graded components of the canonical module are true, but
there is a close relationship between the two conditions. In fact in
one of our main results we express, more generally, the annihilators
of graded components of $\omega$ in terms of colon ideals of powers
of $\m$,
\[{\rm
ann}_R([\omega]_{\leq t}R)={\rm ann}_R([\omega]_{t}R)=\oplus_i
[J^{i+j-d+t+1} : \m^{j}]_i\] for every $t$ and $j \geq a+d$ (see
Theorem~\ref{ann of omega}). Conversely, this allows one to write
the colon ideals $J^i:\m^j$ for any $i\/$ and $j\geq a+d$ in the
form
\[ J^{i} : \m^{j} =\m^{i-j+a+d} + N \, , \]
where $\m^{i-j+a+d}$ is the `expected' part and $N$ is an ideal of
height zero that is generated in degrees $\leq {i-j+a+d-1}$ and can
be described as
\[N= (\bigoplus_{l \leq i-j+a+d-1} [{\rm ann}_R
([\omega]_{i-j-l+d-1} R)]_{l})R \,
 \]
 (see Corollaries \ref{colondescription} and \ref{K}). To prove
 these results it is useful to replace the colon ideals $J^{i} :_R
 \m^{j}$
 in $R$ by the corresponding colon ideals $J^{i}\omega :_{\omega} \m^{j}$ in $\omega$,
 which we then relate to truncations $[\omega]_{\geq i-j+d}$ of
 $\omega$ and to graded components of the canonical module $\Omega$
 of the extended Rees ring of $\m$.
This is done in one of our main technical results. There we also
use a bound on the regularity of $\Omega$ to show that
$J^{i}\omega :_{\omega} \m^{j}= J^{i-j}(J^{j}\omega
:_{\omega} \m^{j})$ for {\it every} integer $i$ with $i \geq j$
(see Theorem~\ref{omega} and, for related results, \cite{HySm1, PU,
HySm2, CC}).

Our results imply, for instance, that $[\omega]_{-a}R$ is faithful,
equivalently, $[\omega]_{t}R$ is faithful for every $t \geq -a$, if and only
if $J^{i} : \m^{j} =\m^{i-j+a+d}$ for every $i$ and $j \geq a+d$, if
and only if $J^{i} : \m^{a+d} =\m^{i}$ for some $i \gg 0$ (see
Corollary~\ref{omegafaithful2}). The question arises how large the
integer $i$ has to be chosen in the last statement. Conceivably
$i=a+d+1$ always works at least when $k$ is perfect and $R$ is reduced,
and we can show this if $d=1$ or if $R$ is an almost complete
intersection of embedding codimension $2$, for instance (see
Corollaries~\ref{omegafaithfuldim1} and \ref{codim2}). In general we
can prove that $i=\alpha +d+1$ suffices, where $\alpha=\alpha(R)$ is
the smallest possible $a$-invariant of a standard graded Gorenstein
ring of dimension $d$ mapping onto $R$ so that the kernel vanishes locally at
its minimal primes (see Corollary~\ref{alpha2}).
We deduce this result from an estimate on the initial degree of
certain colon ideals (see Proposition~\ref{primes}), which in turn
follows from a general bound on the generic generator degree of the canonical module (see
Proposition~\ref{CU}). The stronger result in the case of almost complete
intersections of embedding codimension $2$ is proved using an estimate for
the generic generator degree of the first syzygy module of homogeneous ideals
(see Propositions~\ref{syzygies}).

Returning to the core and the general assumptions used in this
context, we conclude that $[\omega]_{-a}R\/$ is faithful if and only
if $\core{\m^n}=\m^{nd+a+1}$ for every $n\geq 1$, if and only if
$\core{\m^n}$ is generated in one degree for some $n \gg 0$ (see
Theorem~\ref{core-ann}). Furthermore, $\core{\m^n}$ can be replaced
by $\core{\m}$ in the last statement, provided $d=1$ or $R$ is a
reduced almost complete intersection of embedding codimension $2$,
for instance (see Corollary~\ref{indeg=a+d}). This leads, rather
directly, to a geometric interpretation of the core
when $R$ is the homogeneous coordinate ring of a finite set $X$ of
reduced points in projective space. Most notably,
$\core{\m}=\m^{a+2}$ if and only if $X$ has the Cayley-Bacharach
property (see Corollary~\ref{CB}). Recall that $X$ is {\it
Cayley-Bacharach} if the Hilbert function of $X\setminus \{P\}$ does
not depend on the point $P \in X$. Since this
property is equivalent to the faithfulness of $[\omega]_{-a}R$ (see
\cite{GKR}), the above characterization in terms of the core then
follows as an immediate consequence of Corollary~\ref{indeg=a+d}. We
also show that if a large enough subset of $X$ lies on a
hypersurface of low degree then the initial degree of $\core \m$ is
forced to be unexpectedly small (see Corollary~\ref{Y and Z} and
Proposition~\ref{local}), underlining once more the fact that
the shape of the core reflects uniformity properties of the set of points.
%are reflected in the `expected form' of $\core{\m}$.

\s

\section{What annihilates the components of the canonical module?}

\sss

We begin by fixing notation and recalling some general facts.
Let $k$ be a field and $R$ a standard graded
$k$-algebra of dimension $d$ with homogeneous maximal ideal
$\m$ and graded canonical module $\omega$. Recall that
\begin{align*}
 \hspace{-1.5cm} a(R)&= - {\rm
 min} \, \{ \, i \, | \,  [\omega]_{i} \not=0 \}
 \end{align*}
is the $a$-{\it invariant}  of $R$. We also consider the integers
\begin{align*}
b(R)&=-{\rm min} \{ i \,  | \, [\omega]_iR \, \mbox{\rm {is
$R$-faithful}}
\} \ \  {\rm and} \\
c(R)&=-{\rm max} \{ i \,  | \, [k \otimes_{R} \omega]_i \not=0 \,
\}.
\end{align*}
By local duality, the $a$-invariant is the top degree of the local
cohomology module $H^d_{\m}(R)$. Also notice that $c(R)$ is the
negative of the largest generator degree of $\omega$. Since $\omega$
is $R$-faithful we have $a(R) \ge b(R) \ge c(R)$.

Now assume that $R$ is Cohen-Macaulay.
Let $S=k[x_1, \ldots, x_n]$ be a polynomial ring mapping
homogeneously onto $R$ and write $g={\rm codim}_S \, R = {\rm
projdim}_S \, R$. The integers $a(R)$ and $c(R)$ can be expressed in
terms of the minimal homogeneous free $S$-resolution of $R$ by means
of the formulas
\begin{align*}
a(R)&= {\rm max} \, \{ \, i \, | \, [{\rm Tor}^S_g(k,R)]_i \not=0 \,\} -n  \\
c(R)&= \, {\rm min} \, \{ \, i \, | \, [{\rm Tor}^S_g(k,R)]_i \not=0
\,\} -n \, .
\end{align*}
\noindent
Notice that $c(R) \geq -d$; in particular, $\omega$ is
generated in degrees at most $d$ and $a(R) + d \geq 0$.

We also consider the extended Rees algebra $R[\m t,t^{-1}]$, which
is a bigraded subring of $R[t,t^{-1}]$. There are natural maps
of bigraded modules
\[
\omega_{R[\m  t, t^{-1}]} \hookrightarrow (\omega_{R[\m  t,
t^{-1}]})_{t^{-1}} \cong \omega_{R[\m t, t^{-1}]_{t^{-1}}} =
\omega_{R[t, t^{-1}]} \cong \omega \otimes_{R} R[t,t^{-1}]= \oplus\,
\omega \, t^{i}.
\]
Identifying $\omega_{R[\m t, t^{-1}]}$ with its image in
$\oplus\, \omega \, t^{i}$ we obtain a bigraded canonical module
\[
\Omega= \oplus\, \Omega_i t^{i} \subset \oplus\, \omega \, t^{i},
\]
so that $\Omega_i=\omega$ for $i \ll 0$. One has homogenous
isomorphisms $ R[\m \/ t, t^{-1}]/(t^{-1}) \cong {\rm gr}_{\m} (R)
\cong R$, thinking of the last $R$ as being diagonally bigraded.
They induce identifications
\begin{equation}\label{omegarees}
  \Omega/t^{-1}\Omega \cong \omega (1) \ \ \hbox{\rm and} \ \
\Omega_i/\Omega_{i+1} \cong [\omega]_{i+1}.
\end{equation}

We will use the convention that the power of any element or ideal with
non-positive exponent is one or the unit ideal, respectively.

%ssume in addition that $d\geq 1$. Let $y_1, \ldots, y_d$
%e a system of parameters in $R$ consisting of
%inear forms, which always exists if $k$ is infinite.
%he ideal $J=(y_1, \ldots, y_d)$
%s a minimal reduction of $\m$ with reduction number $r_J(\m)=a(R)+d$, as can be seen by reducing modulo $J$.
%herefore
%[ (y_1, \ldots , y_d)^{i}:\m^{j}= (y_1^{i}, \ldots, y_d^{i}):\m^{j+(i-1)(d-1)} \] for every $i$ and $j \geq a(R)+d$;
%n fact, one can show as in \cite[proof of 2.2]{PUV} that if $I$ is any
%\m$-primary ideal and $\beta_1, \ldots , \beta_d$ are elements of $I$
%ith $I^{j+1}=(\beta_1, \ldots , \beta_d)I^j$ for some integer $j$, then
%begin{equation}\label{PUV 2.2}(\beta_1,\ldots,\beta_d)^{i}:I^{j}=(\beta_1^{i},\ldots,\beta_d^{i}): I^{j+(i-1)(d-1)}
%end{equation}
%or every $i$.

\ms

\begin{Assumptions}\label{standard}{\rm We assume $k$ is a field and $R$ is
a standard graded Cohen-Macaulay $k$-algebra of dimension $d\geq 1$ with homogeneous
maximal ideal $\m$ and graded canonical module $\omega$. We write $a=a(R)$, $b=b(R)$, $c=c(R)$,
and let $\Omega = \omega_{R[\m t,t^{-1}]}$ be as above, with $\Omega_i=\omega$ for $i \ll 0$.
Let $y_1, \ldots, y_d$ be a system of parameters in $R$ consisting of
linear forms. We write $J$ for the ideal generated by $y_1, \ldots, y_d$ and $J^{[i]}$ for the ideal
generated by the powers $y_1^{i}, \ldots, y_d^{i}$, where $i$ is any integer.}
\end{Assumptions}

\s

Elements $y_1, \ldots, y_d$ as in Assumptions \ref{standard} always
exist if $k$ is infinite. The ideal $J$ they generate is a minimal
reduction of $\m$ with reduction number $r_J(\m)=a+d$, as can be
seen by reducing modulo $J$. Therefore
\[ J^{i}:\m^{j}= J^{[i]}:\m^{j+(i-1)(d-1)} \] for every $i$ and $j \geq a +d$;
in fact, one can show as in \cite[proof of 2.2]{PUV} that if $I$ is any
$\m$-primary ideal and $\beta_1, \ldots , \beta_d$ are elements of $I$
with $I^{j+1}=(\beta_1, \ldots , \beta_d)I^j$ for some integer $j$, then
\begin{equation}\label{PUV 2.2}(\beta_1,\ldots,\beta_d)^{i}:I^{j}=(\beta_1^{i},\ldots,\beta_d^{i}): I^{j+(i-1)(d-1)}
\end{equation}
for every $i$.

\ms

In this section we prove our results about the relationship between
annihilators of graded components of $\omega$ on the one hand and
the colon ideals $J^i : \m^j$ on the other hand. Here it suffices to
consider $j=a+d$, as the next remark shows.

\ms

\begin{remark}\label{colon1}
With assumptions as in \ref{standard} one has
\[ J^{i} : \m^{j}=J^{i-j+a+d}:\m^{a+d} \] for every $i$ and $j
\geq a+d$.
\end{remark}
\proof Since $j \geq a+d=r_J(\m)$, we obtain
\[
J^{i}:\m^{j}=J^{i}:J^{j-a-d}\m^{a+d}=(J^{i}:J^{j-a-d}):
\m^{a+d}=J^{i-j+ a+d}:\m^{a+d}.
\]
The last equality holds because the associated graded ring of $J$ has positive depth.
\QED

\s

We first study the colons $J^i \omega :_{\omega} \m^j$ in $\omega$,
which exhibit a more regular behavior than the corresponding ideals
$J^i :_R \m^j$.

\begin{theorem}\label{omega}
In addition to the assumptions of $\  \ref{standard}$ let $i$ and $j
\geq a+d$ be integers, and set $s=j+(i-1)(d-1)$. One has
\begin{itemize}
  \item[(a)] $J^{i}\omega :_{\omega}
  \m^{j} =J^{[i]}\omega :_{\omega} \m^s = [\omega]_{\geq i-j+d} =
  \Omega_{i-j+d-1}$
  \item[(b)]$J^{i}\omega :_{\omega}
  \m^{j}= J^{i-j}(J^{j}\omega :_{\omega}
  \m^{j})$ whenever $i \geq j$.
\end{itemize}
\end{theorem}
\proof To prove part (a) we first notice that
\[ \Omega_{i-j+d-1}=J^{i}\omega :_{\omega}
  \m^{j}
\]
according to \cite[3.7]{CC}. Next,
\begin{eqnarray*}
\nonumber %to remove numbering (before each equation)
 J^{i}\omega :_{\omega} \m^{j} &=& (J^{[i]}:_R J^{(i-1)(d-1)})
 \omega :_{\omega} \m^j  \hspace{1.55cm} \hbox{\rm by (\ref{PUV 2.2})} \\
   &\subset&  (J^{[i]}\omega :_{\omega} J^{(i-1)(d-1)}):_{\omega} \m^j \\
   &=&  J^{[i]}\omega :_{\omega} (J^{(i-1)(d-1)} \m^j ) \\
   &=& J^{[i]}\omega :_{\omega} \m^s   \hspace{4cm} \hbox{\rm because $j \geq a+d = r_{J} (\m)$} \, .\\ \end{eqnarray*}

To show that $J^{[i]}\omega :_{\omega} \m^s= [\omega]_{\geq i-j+d}$
we may reduce modulo  $J^{[i]}\omega\,$; indeed, $J^{[i]}\omega \subset
[\omega]_{\geq i-j+d}$ because $\omega$ is concentrated in degrees
$\geq -a$ and $i-a \geq i-j+d$. However, $(J^{[i]}\omega :_{\omega}
\m^s)/ J^{[i]}\omega \cong \ol{0} :_{\ol{\omega}} \ol{\m}^s$, where
$\ol{\m}$ is the homogeneous maximal ideal of $\,  \ol{R}=R/J^{[i]}$ and
$\ol{\omega}= \omega/J^{[i]}\omega \cong \omega_{\ol{R}} \,(-id)$.
As $\omega_{\ol{R}}$ is an Artinian module with socle concentrated
in degree $0$, the module $\ol{\omega}$ is Artinian with socle
concentrated in degree $id$. Thus $\ol{0} :_{\ol{\omega}} \ol{\m}^s
= [\ol{\omega}]_{\geq id+1-s}= [\ol{\omega}]_{\geq i-j+d}$. It
follows that $J^{[i]}\omega :_{\omega} \m^s= [\omega]_{\geq i-j+d}$.

So far we have shown the inclusions
\[
\Omega_{i-j+d-1}=J^{i}\omega :_{\omega}
  \m^{j} \subset J^{[i]}\omega :_{\omega} \m^s = [\omega]_{\geq
  i-j+d} \, .
\]
Now $\Omega =\oplus \Omega_{l}t^l \subset \oplus  [\omega]_{\geq
l+1} t^{l}=\Omega'$ are bigraded $R[\m t, t^{-1}]$-modules that are
finitely generated because $[\omega]_{\geq l+1}=\m [\omega]_{\geq
l}$ for $l \gg 0$. By (\ref{omegarees}) this inclusion induces an
isomorphism $\Omega/t^{-1}\Omega \cong  \Omega'/t^{-1}\Omega'$,
which shows that $\Omega'=\Omega + t^{-1}\Omega'$. Therefore
$\Omega'=\Omega$ by the graded Nakayama lemma.

We now prove part (b). In the light of (a) we need to show that for
every $l \geq d$ one has $[\omega]_{\geq l+1} = J [\omega]_{\geq
l}$, or equivalently, $[\omega]_{\geq l+1} \subset J \omega$.
However, $\omega/ J \omega \cong \omega_{R/J} (-d)$ is concentrated
in degrees $\leq d$, which gives $[\omega]_{\geq l+1} \subset J
\omega$. \QED

\s

The next result addresses the comparison between the colons
$J^i\omega :_{\omega} \m^j$ and $J^i :_R \m^j$.

\begin{corollary}\label{=}
In addition to the assumptions of $\,  \ref{standard}$ let $i$ and $j
\geq a+d$ be integers. If \[J^{i}\omega :_{\omega} \m^{j}=(J^{i}:_R
\m^{j})\, \omega \hspace{.5cm} \hbox{\rm for some} \hspace{.3cm} i \geq
j \, ,\]  then for every $l \geq i$,
\begin{itemize}
  \item[(a)] $ J^{l}\omega :_{\omega} \m^{j}=(J^{l}:_R \m^{j}) \, \omega$
  \item[(b)] $J^{l}:_R \m^{j}$ is integral over  $J^{l-i}(J^{i}:_R
  \m^{j})$; in particular, if $R$ is reduced then the initial degree
  of the former ideal satisfies \[{\rm indeg} (J^{l}:_R \m^{j}) ={\rm indeg} (J^{i}:_R
  \m^{j})+ l-i \, .\]
\end{itemize}
\end{corollary}
\proof One has
\begin{eqnarray*}
 \nonumber %to remove numbering (before each equation)
  (J^{l}:_R \m^j) \, \omega &\subset& J^{l}\omega :_{\omega} \m^j\\
   &=& J^{l-i}(J^{i}\omega :_{\omega} \m^j) \hspace{2cm} \hbox{\rm by Theorem~\ref{omega}(b)} \\
  &=& J^{l-i} (J^{i}:_R \m^j) \, \omega \hspace{1.935cm} \hbox{\rm by our assumption}\\
 &\subset& (J^{l}:_R \m^j) \, \omega \, .
\end{eqnarray*}
This immediately gives $ J^{l}\omega :_{\omega} \m^{j}=(J^{l}:_R
\m^{j}) \, \omega$, proving assertion (a). Furthermore, $(J^{l}:_R
\m^j) \, \omega=J^{l-i} (J^{i}:_R \m^j) \, \omega$, which implies (b) since
$\omega$ is a faithful $R$-module. \QED

\s

The next proposition gives a characterization for when the equality
assumed in the previous corollary obtains.

\begin{proposition}\label{omegacontainment}
In addition to the assumptions of $\, \ref{standard}$ let $i$ and $j \geq
a+d$ be integers. Set $s=j+(i-1)(d-1)$ and let $^{^{\tratto}}$ denote images
in $\ol{R}=R/J^{[i]}$. One has
\[J^{i}\omega :_{\omega}
\m^{j}=(J^{i} :_{R} \m^{j}) \, \omega \]
if and only if
\[\ol{\m}^{s}=
{\ol{0}:_{\ol{R}}(\ol{0}:_{\ol{R}}\ol{{\f{m}}}^{s})} \, .\]

\end{proposition}
\proof First notice that $\ol{R}$ is an Artinian ring
and that $\ol{\omega}= \omega/J^{[i]}\omega \cong
\omega_{\ol{R}} \,(-id)$. From Theorem~\ref{omega}(a)
one has $J^i\omega :_{\omega} \m^j = J^{[i]}\omega :_{\omega}
\m^s$. Hence it follows that $(J^i\omega :_{\omega} \m^j)/
J^{[i]}\omega \cong \ol{0} :_{\ol{\omega}}
\ol{\m}^s$. Now the latter module is naturally isomorphic to
$(\ol{0} :_{\omega _{\ol R}} \ol{\m}^s)(-id) \cong \omega_{\ol{R}/\ol{\m}^s} \,(-id)$.
Thus we have shown that
\[ (J^i\omega :_{\omega} \m^j)/
J^{[i]}\omega \cong \omega_{\ol{R}/\ol{\m}^s} \,(-id).\]

On the other hand, from (\ref{PUV 2.2}) one knows that $J^{i} : \m^{j}=
J^{[i]}:\m^{s}$. Therefore $(J^{i} :_R \m^{j}) \, \omega/
J^{[i]}\omega =(\ol{0} :_{\ol{R}} \ol{\m}^s) \, \ol{\omega}$, which is
isomorphic to
$(\ol{0} :_{\ol{R}} \ol{\m}^s)\, \omega_{\ol{R}} \,(-id)$. The last module
can be identified with
$\omega_{\ol{R}/{\ol{0}:(\ol{0}:\ol{\m}^{s})}}(-id)$, as shown in \cite[2.3(b)]{U94}.
Hence we have proved that
\[(J^{i} :_R \m^{j}) \, \omega/ J^{[i]}\omega \cong \omega_{\ol{R}/{\ol{0}:(\ol{0}:\ol{\m}^{s})}}(-id) \, . \]

We conclude that the obvious containment $(J^{i}:_R \m^{j}) \, \omega
\subset J^{i}\omega :_{\omega} \m^{j}$ is an equality if and only if
the inclusion $\ol{\m}^{s} \subset \ol{0}:(\ol{0}:\ol{\m}^{s})$ is.
\QED

\s

If the equalities of Proposition~\ref{omegacontainment} hold for
all $i$, then the canonical module of the extended Rees ring of $\m$ has
an easy description, namely
\[\Omega=(R[Jt,t^{-1}] :_{R[t,t^{-1}]} \m^j) t^{d-1-j}\omega\]
(see Theorem ~\ref{omega}(a) or \cite[3.7]{CC}).
These equalities obtain, quite generally, when $i \geq j$ and $d=1$:

\begin{proposition}\label{=indim1}
In addition to the assumptions of $\, \ref{standard}$ suppose that $d=1$
and $R$ is reduced. For
every $i \geq j$ and $j \geq a+1$ one
has \[(J^{i} :_R \m^{j}) \, \omega =  J^{i}\omega :_{\omega}
\m^{j}=[\omega]_{\geq i-j+1}.\]
\end{proposition}
\smallskip
\proof The second equality follows from Theorem~\ref{omega}(a). We
prove the first equality. It suffices to consider the case $i=j=a+1$
according to Remark \ref{colon1} and Corollary~\ref{=}(a),
or simply because $J$ is generated by a single regular element.
In light of Proposition~\ref{omegacontainment} we need to show that
\[\m^{a+1}=y^{a+1}R:_R ( y^{a+1}R :_R \m^{a+1})\]
with $y=y_1$.

In the total ring of quotients $K$ of $R$ we consider the
%blowup
%$B=R[\m/y]$ of $\m$. Write $e$ for the multiplicity of $R$. Passing
%to the Artinian ring $R/yR$ one sees that the Hilbert function
%${\rm{H}}_R(n)$ attains its maximal value $e$ if and only if $n \geq
%a+1$. Since $a+1=r_J(\m)$ one has
%$B = \m^{a+1}/y^{a+1} \cong \m^{a+1}(a+1) = R_{\geq a+1}(a+)$. Thus $B$ is a
%non-negatively graded finite $R$-module with constant Hilbert
%function, ${\rm{H}}_B(n)=e$ for $n\geq 0$. It follows that the
%conductor $R:_K B \/$ is $B_{\geq a+1} = R_{\geq a+1} =\m^{a+1}$.
%Hence we obtain
%\begin{eqnarray*}
%\nonumber
%y^{a+1}R:_R ( y^{a+1}R :_R \m^{a+1}) &=& y^{a+1}R:_K ( y^{a+1}R :_K \m^{a+1}) \\
%   &=& R:_K (R :_K \m^{a+1})  \\
%   &=& R:_K ( R :_K (R:_KB))      \hspace{0.95cm}  \hbox{\rm since} \ R:_KB = \m^{a+1} \\
%   &=& R:_K B \\
%   &=& \m^{a+1}  \hspace{3.5cm} \hbox{\rm since} \ R:_KB = \m^{a+1} .
integral closure $S$ of $R$. Our assumptions, most
notably the reducedness of $R$, imply that $R \subset S$ is a
homogeneous inclusion of non-negatively graded Noetherian rings and
moreover $R:_K ( R :_K S)=S$. Thus we obtain
\begin{eqnarray*}
\nonumber
 \m^{a+1} &\subset& y^{a+1}R:_R ( y^{a+1}R :_R \m^{a+1}) \\
   &\subset& y^{a+1}R:_R ( y^{a+1}R :_R \m^{a+1}S) \\
   &=& y^{a+1}R:_R ( y^{a+1}R :_R y^{a+1}S)  \hspace{1.5cm}\hbox{\rm
since $yR$ is a reduction of $\m$}\\
   &=& y^{a+1}R:_R ( R :_R S)      \hspace{2.925cm}  \hbox{\rm since $y$
is $R$-regular} \\
   &=& [y^{a+1}R:_K ( R :_R S)] \cap R \\
   &=& [y^{a+1}(R:_K ( R :_K S))] \cap R \\
   &=& y^{a+1}S \cap R     \hspace{4.17cm} \hbox{\rm since $R:_K ( R :_K
S)=S$}\\
   &\subset&  \m^{a+1}  \hspace{4.97cm} \hbox{\rm since $S$ is non-negatively
   graded}.
\end{eqnarray*}\QED

\s

Next we are going to use the above results to say something about
the annihilators of the graded components of $\omega$.

\begin{theorem}\label{ann of omega}
With assumptions as in $\ref{standard}$ one has \[{\rm
ann}_R([\omega]_{\leq t}R)={\rm ann}_R([\omega]_{t}R)=\oplus_i
[J^{i+j-d+t+1} : \m^{j}]_i\] for every $t$ and $j \geq a+d$.
\end{theorem}
\proof The first equality is obvious because $R$ contains a linear
form, namely $y_1$, that is regular on $\omega$.

Theorem~\ref{omega}(a) shows that
\[(J^{i+j-d+t+1} : \m^{j}) \, \omega \subset [\omega]_{\geq i+t+1} \, ,\]
which gives $[J^{i+j-d+t+1} : \m^{j}]_i [\omega]_{\leq t} =0.$ This
proves the inclusion
\[\oplus_i [J^{i+j-d+t+1} : \m^{j}]_i \subset {\rm
ann}_R([\omega]_{\leq t}R) \, .\]

To show the containment \[{\rm ann}_R([\omega]_{\leq t}R) \subset
\oplus_i [J^{i+j-d+t+1} : \m^{j}]_i \, ,\] we choose an element $f\in
[{\rm ann}_R([\omega]_{\leq t}R)]_i.$ We need to prove that \[ f \in
J^{i+j-d+t+1} : \m^{j} = J^{[i+j-d+t+1]} : \m^{j+ (i+j-d+t)(d-1)} \, ,\]
where the last equality holds by (\ref{PUV 2.2}). Write $s=j+
(i+j-d+t)(d-1)$ and let $^{^{\tratto}}$ denote images in the Artinian ring
$\ol{R}=R/J^{[i+j-d+t+1]}.$ The asserted inclusion $f \m^s \subset
J^{[i+j-d+t+1]}$ is equivalent to $\ol{f}\, \ol{\f{m}}^{s} = \ol{0},$
which in turn means that $\ol{f} \,\ol{\f{m}}^{s} \omega_{\ol{R}} =
\ol{0}$ as $\omega_{\ol{R}}$ is faithful over $\ol{R}$. Since
$\omega_{\ol{R}} = \ol{R} \otimes_R  \omega \, ((i+j-d+t+1)d)$ and
$f [\omega]_{\leq t}=0$ it follows that $\ol{f}
[\omega_{\ol{R}}]_{\leq t-(i+j-d+t+1)d}=\ol{0}.$ Therefore $\ol{f}
[\omega_{\ol{R}}]_{\leq -s-i}=\ol{0},$ and as $\ol{f}$ has degree
$i$ we conclude that \[ \ol{f} \,\ol{\m}^s \omega_{\ol{R}} \subset
[\omega_{\ol{R}}]_{\geq 1} =\ol{0} \, . \]\QED

\s

Again, one has a better result in dimension one.

\begin{corollary}\label{ann of omega dim1}
In addition to the assumptions of $\, \ref{standard}$
suppose that $d=1$ and write $y=y_1$. One has
\[
{\rm ann}_R([\omega]_{\leq t}R)={\rm ann}_R([\omega]_{t}R)=\oplus_{i
\leq -t-1} [J^{i+j+t} : \m^{j}]_i \, \bigoplus \ [J^j:\m^j]_{-t}k[y]
\] for every $t$ and $j \geq a+1$.
\end{corollary}
\proof We apply Theorem~\ref{ann of omega} using the fact that
$J^{i+j+t} : \m^j= y^{i+t}(J^j:\m^j)$ for every $i \geq -t$. \QED

\s

The above corollary shows in particular that if $d=1$ then the
graded $R$-module ${\rm ann}_R([\omega]_{\leq t}R)={\rm ann}_R([\omega]_{t}R)$
has Castelnuovo-Mumford regularity at most $-t$; see also Proposition
\ref{CM}. For the definition and basic properties of the Castelnuovo-Mumford
regularity we refer to \cite[p.168]{BH}.
% and \cite[Chapter 20.5]{ei}.

\sss

We are now ready to answer one the main questions raised in the
introduction. With the next three corollaries we characterize the
faithfulness of submodules generated by graded components of
$\omega$, in terms of certain colon ideals of powers of $\m$ having
an `expected form'. In light of Remark \ref{colon1} we may restrict
ourselves to the case $j=a+d$.

\ms
 \begin{corollary}\label{omegafaithful1} With assumptions as in $\ref{standard}$ the
 following are equivalent for an integer $t :$
 \begin{itemize}
 \item[(a)] $[\omega]_{t}R$ is faithful, i.e., $t \geq -b$
 \item [(b)]$J^{i+a+t} :
\m^{a+d} \subset \m^{i}$  for every  $ i$
\item[(c)]$J^{i+a+t} :
\m^{a+d} \subset \m^{i}$  for some $i \gg 0$.
\end{itemize}
If $R$ is reduced and \, $J^{i+a+t}\omega :_{\omega}
\m^{a+d}= (J^{i+a+t} :_R
\m^{a+d}) \, \omega$ \, for some $i \geq d-t$, then the above conditions
are equivalent to $:$
\begin{itemize}
\item[(d)] $J^{i+a+t} : \m^{a+d} \subset \m^{i}$.
\end{itemize}
\end{corollary}
\proof Recall that $R$ contains a linear form that
is a regular element.  Thus, if a homogeneous ideal vanishes in
a certain degree it also vanishes in every smaller degree. Now
Theorem~\ref{ann of omega} gives the equivalence of (a) and
(b). The same theorem shows that if (c) holds then $[{\rm ann}_R([\omega]_{t}R)]_{i-1}=0$ for some
$i-1 \gg 0$, proving (a). Finally, (d) implies (c) according to Corollary~\ref{=}(b). \QED

\ms

 \begin{corollary}\label{omegafaithful2} With assumptions as in $\ref{standard}$ the
 following are equivalent $:$
 \begin{itemize}
 \item[(a)] $[\omega]_{-a}R$ is faithful
 \item [(b)]$J^{i} :
\m^{a+d}=\m^{i}$  for every  $ i$
\item[(c)]$J^{i} :
\m^{a+d}=\m^{i}$  for some $i \gg 0$.
\end{itemize}
If $R$ is reduced and  $J^{i}\omega :_{\omega}
\m^{a+d} =(J^{i} :_R
\m^{a+d}) \, \omega $ for some $i \geq a+d$, then the above conditions
are equivalent to $:$
\begin{itemize}
\item[(d)] $J^{i} : \m^{a+d} = \m^{i}$.
\end{itemize}
\end{corollary}
\proof Notice that $\m^{i} \subset J^{i}:\m^{a+d}$ for every $i\/$
because $a+d =r_J(\m)$. Now the assertions follow from
Corollary~\ref{omegafaithful1}.  \QED

\ms

 \begin{corollary}\label{omegafaithfuldim1} In addition to the
 assumptions of $\, \ref{standard}$
 suppose that $d=1$. The
 following are equivalent $:$
 \begin{itemize}
 \item[(a)] $[\omega]_{-a}R$ is faithful
 \item [(b)]$J^{i} :
\m^{a+1}=\m^{i}$  for every  $ i$
\item[(c)]$J^{i} :
\m^{a+1}=\m^{i}$  for some $i \geq a + 1$.
\end{itemize}
\end{corollary}
\proof In light of Corollary~\ref{omegafaithful2} it suffices to
prove that if  $J^{i} : \m^{a+1}=\m^{i}$  for some $i \geq a +1$,
then $J^{l} : \m^{a+1}=\m^{l}$  for every $l \geq i$. Indeed, since
$J$ is generated by a single regular element and $i\geq a+1= r_J(\m)$, it
follows that \[J^{l} : \m^{a+1}= J^{l-i}(J^{i} : \m^{a+1})=J^{l-i} \m^{i}=\m^{l}. \]\QED

\ms

The faithfulness of $[\omega]_{-a}R$, together with the additional
condition $J^{i}\omega :_{\omega}\m^{a+d} =(J^{i} :_R
\m^{a+d}) \, \omega $ in Corollary \ref{omegafaithful2}, means that
$\omega$ and $[\omega]_{-a}R$ coincide from degree $i-a$ on:

\sss

\begin{remark}\label{==}{\rm In addition to the assumptions of \ref{standard} let
$i \geq a+d$ be an integer. Then the equality $ [\omega]_{\geq i-a} =
\m^i [\omega]_{-a}$ holds if and only if $J^{i}\omega :_{\omega}\m^{a+d}
=(J^{i} :_R \m^{a+d}) \, \omega $ and $[\omega]_{-a}R$ is faithful.
The forward direction follows because
$[\omega]_{\geq i-a} =J^{i}\omega :_{\omega}\m^{a+d}$ according to
Theorem~\ref{omega}(a), $\m^i \subset J^i : \m^{a+d}$, and $[\omega]_{\geq i-a}$ is
faithful. For the converse notice that
%the faithfulness of $[\omega]_{-a}R$ and
Corollary~\ref{omegafaithful2} gives $J^i:\m^{a+d}=\m^i$, hence $[\omega]_{\geq i-a}=
J^{i}\omega :_{\omega}\m^{a+d} =\m^i\omega$.
Therefore $[\omega]_{\geq i-a}=\m^i
[\omega]_{-a} + [\omega]_{\geq i+1-a}.$ However, $\omega$ is generated in
degrees at most $d \leq i-a$, hence $[\omega]_{\geq
i-a}$ is generated in degree $i-a$. It follows that $[\omega]_{\geq i-a} = \m^i
[\omega]_{-a}$.}
\end{remark}
\ms
Unfortunately, the faithfulness of $[\omega]_{-a}R$ alone is not
sufficient to guarantee the equality $J^{i}\omega :_{\omega}\m^{a+d} =
(J^{i} :_R \m^{a+d}) \, \omega$, as the next example shows.
\sss

\begin{example} In addition to the assumptions of \ref{standard} suppose
that $R$ is a domain of type 2 which is not level and which is not
Gorenstein locally on the punctured spectrum. Clearly
$[\omega]_{-a}R$ is faithful since $R$ is a domain. However,
$J^{i}\omega :_{\omega} \m^{a+d} \not= (J^{i} :_R \m^{a+d}) \, \omega$ for every
$i\geq a+d$, because otherwise Remark~\ref{==} implies that $\omega$
and the cyclic module $[\omega]_{-a}R$ would coincide locally on the
punctured spectrum.
\end{example}

\s

Conversely, Theorem~\ref{ann of omega} can be used to obtain
information about the colon ideals $J^i : \m^j\,$. This is done in
the remaining corollaries of this section.

\ms

 \begin{corollary}\label{colondescription}
With assumptions as in $\ref{standard}$ one has
\[ J^{i} : \m^{j} =\m^{i-j+a+d} \ \bigoplus \ \ \ \bigoplus_{l =
i-j+b+d}^{i-j+a+d-1} [{\rm ann}_R ([\omega]_{i-j-l+d-1} R)]_{l}
 \]
for every $i$ and $j\geq a+d$.
\end{corollary}
\proof Theorem~\ref{ann of omega} shows that
\[ [{\rm ann}_R ([\omega]_{i-j-l+d-1} R)]_{l}= [J^{i} : \m^{j}]_{l}
\] for every $l$. On the other hand,
$[{\rm ann}_R ([\omega]_{i-j-l+d-1} R)]_{l}=0 \/$ for $l\leq
i-j+b+d-1$ by the definition of $b$. Finally, $\m^{i-j+a+d} \subset
J^i : \m^j$. \QED

\ms

\begin{corollary}\label{K} With assumptions as in $\ref{standard}$ there exists
an ideal $N$ of height zero such that
\[ J^{i} : \m^{j} =\m^{i-j+a+d} + N
 \]
for every $i$ and $j\geq a+d$.
\end{corollary}
\proof Set $N=(\bigoplus_{l = i-j+b+d}^{i-j+a+d-1} [{\rm ann}_R
([\omega]_{i-j-l+d-1} R)]_{l})R.\,$ According to
Corollary~\ref{colondescription}, $J^i: \m^j = \m^{i-j+a+d} + N$. On
the other hand, the first equality of Theorem~\ref{ann of omega}
shows that $N \subset {\rm ann}_R ([\omega]_{-a} R ). \/$ The latter
ideal has height zero since $[\omega]_{-a} R$ is a nonzero submodule
of the maximal Cohen-Macaulay module $\omega$. \QED

\ms

\begin{corollary}\label{power} In addition to the assumptions of $\, \ref{standard}$
let $i$ and $j \geq a+d$ be integers. If $J^i : \m^j$ is generated
in one degree then $J^i : \m^j=\m^{i-j+a+d}$.
\end{corollary}
\proof The claim follows from Corollary~\ref{K} since $J^i : \m^j$ has positive
height. \QED

\ms

\begin{corollary}\label{colonmax1} With assumptions as in $\ref{standard}$ one has
\[ \m^{i-j+a+d}  \subset J^{i} : \m^{j} \subset \m^{i-j+b+d}.
\]
for every $i$ and $j\geq a+d$.
\end{corollary}
\proof The containments are direct consequences of
Corollary~\ref{colondescription}. \QED

\ms

\begin{corollary}\label{colonmax2}
In addition to the assumptions of $\, \ref{standard}$ assume that $R$ is
a domain or $R$ is level. One has
\[
J^{i} \colon \m^j =\m^{i-j+a+d}
\]
for every $i$ and $j \ge a+d$.
\end{corollary}
\proof The result follows from Corollaries~\ref{K} and
\ref{colonmax1}.\QED

\s

In the one-dimensional case the colon ideals $J^i : \m^j$ can be
expressed in terms of conductor ideals:

\begin{remark}\label{coreandS}
{\rm In addition to the assumptions of \ref{standard} suppose that
$d=1$, write $y=y_1$, let $K$ denote the total ring of quotients of
$R$, and let $R[\m/y] \subset B$ be a homogeneous inclusion where $B$ is a non-negatively graded finite
$R[\m/y]$-module contained in $K$. One has $B=R[\m/y]$ and
\[J^i :_R \m^j= J^{i-j}(R :_K B)= \m^{i-j}(R :_K B) \] for every $i\geq j$ and $j\geq a+1$.

To prove the first claim notice that $B$ is a finite $R$-module and hence $B_{\geq l}=R_{\geq l}$
for $l \gg 0$.  Therefore $B_{\geq l}=\m^{l}$, which gives
$B \subset \m^{l}/y^{l}$. As $ j \geq a+1 =r_J(\m)$
one has $\m^{l}/y^{l} = \m^{j}/y^{j}$.
It follows that $B= \m^{j}/y^{j} = R[\m/y]$. The remaining assertions
obtain because $J^j:_R\m^j=y^jR:_K\m^j=R:_KB$ and this is a $B$-module (see also \cite[proof of 3.2]{PU}).}
\end{remark}

\s

\section{Initial degrees of annihilators}

\sss

The previous section leaves open the question of how large the
integer $i$ has to be chosen in Corollaries~\ref{omegafaithful1}(c)
and ~\ref{omegafaithful2}(c). Addressing this issue will entail
establishing an upper bound for the generic generator degree of
first syzygy modules (Proposition~\ref{syzygies}), as well as
estimating the initial degree of certain colon ideals
(Propositions~\ref{primes} and ~\ref{dim1}).

We begin by recalling a result from \cite{CU}, which in turn uses
earlier work on the `fundamental class' (see \cite[p.34]{Elz},
\cite[4.11 and 5.13]{KW}, \cite[3.1]{Li}).

\begin{proposition}\label{CU}\cite[1.1]{CU}
Let $k$ be a perfect field and let $R$ be a standard graded
$k$-algebra of dimension $d$ such that $R_{\p}$ is regular for every
minimal prime $\p$ of dimension $d$. Write $\omega=\omega_R$ and let
$C$ be defined by the exact sequence
\[ 0 \longrightarrow [\omega]_{\leq d}R \stackrel{{\rm nat}}\longrightarrow \omega
\longrightarrow C \longrightarrow 0 \,.
\]
Then ${\rm Supp}(C) \subset {\rm Sing} (R).$
\end{proposition}

\ms

\begin{proposition}\label{syzygies}
Let $k$ be a perfect field and $S$ a standard graded Gorenstein
$k$-algebra of dimension $D$. Let $I$ be an ideal of height $g$ such
that $I_{\P}$ is prime for every prime ideal $\P$ of height $g$ containing
$I$. Assume $I$ is generated by $m \geq g+1$ forms $f_1, \ldots,
f_m$ of degrees $\delta_1 \geq \ldots \geq \delta_m$, and write
$\Delta= \sum_{i=1}^{g+1} \delta_i + D-g +a(S)$. Let $H_{1}$ be the
first Koszul homology of the elements $f_1, \ldots , f_m$ and $C$
the module defined by the exact sequence
\[ 0 \longrightarrow [H_{1}]_{\leq \Delta}S \stackrel{{\rm nat}}\longrightarrow H_{1} \longrightarrow C \longrightarrow 0 \,.
\]
Then $C_{\P}=0$ for every prime $\P$  of $S$ such that
$I_{\P}$ is a complete intersection and $S_{\P}/I_{\P}$ is regular.
If either $g \geq 2$  and $\delta_3 \geq 2$ or $g \geq 1$ and $S$ is
not a polynomial ring, one can replace $H_{1}$ by the first syzygy
module of the elements $f_1, \ldots, f_m$.
\end{proposition}
\proof Fix a prime $\P$ so that $I_{\P}$ is a complete intersection and $S_{\P}/I_{\P}$ is regular.
We may assume that $k$ is infinite and that $f_1, \ldots, f_g$ form
a regular sequence generating $I_{\P}$. Notice
that
\[((f_1, \ldots, f_g) : I)/(f_1, \ldots, f_g) \cong \omega_{S/I}(-\sum_{i=1}^{g} \delta_i
-a(S)).
\]
Therefore Proposition~\ref{CU} shows that there exists a homogeneous element
$\beta$ of degree $D-g +\sum_{i=1}^{g} \delta_i  + a(S)$ in $(f_1,
\ldots, f_g) : I$  generating the factor module $((f_1, \ldots, f_g) :
I)/(f_1, \ldots, f_g)$ locally at  $\P$. In particular,
$\beta \notin \P$ because  $((f_1, \ldots, f_g) : I)_{\P}=S_{\P}$.
Hence for $g+1 \leq j \leq m$ there exist homogeneous syzygies
\[ \beta f_j - \sum_{i=1}^{g} \lambda_{ij} f_i =0
\]
of $f_1, \ldots, f_m$ that have degree at most $\Delta$ and whose
images generate $H_{1}$ locally at $\P$.

Finally, if either $g \geq 2$  and
$\delta_3 \geq 2$ or $g \geq 1$ and $S$ is not a polynomial ring, then the Koszul relations among $f_1,
\ldots, f_m$ have degrees at most $\Delta$.  \QED

\s

We are now going to introduce the $\alpha$-invariant of a graded ring that
will be the basis for many estimates proved in this section.

\begin{definition}\label{DEF}
{\rm Let $k$ be a field and $R$ a standard graded $k$-algebra of
dimension $d$. We define
\[\alpha (R)={\rm min}\{a(S)\} \in {\mathbb{Z}} \cup \, \{\infty\}\/,\]
where $S$ ranges over all standard graded Gorenstein $k$-algebras of
dimension $d$ mapping homogeneously onto $R$ such that
$S_{\P}=R_{\P}$ for every minimal prime $\P \in {\rm Supp}_S(R)$ of
dimension $d$. }\end{definition}

\s
\s
\s

\begin{remark}
\hfill
\begin{itemize}
{\rm \item  [$($a$)$] One has
\[a(R) \leq \alpha(R) \, ,\] because $\omega_R \hookrightarrow \omega_S$
for every $S$ as in Definition~\ref{DEF}.

\item  [$($b$)$] Write $R=k[X_1, \ldots, X_n]/I$ as a factor ring of a
polynomial ring, where $I$ is an ideal of height $g$ generated by
forms of degrees $\delta_1 \geq \ldots \geq \delta_m$. If $\beta_1,
\ldots, \beta_g$ is a homogeneous regular sequence contained in $I$
that generates $I$ at each of its minimal primes of height $g$, then
\[\alpha(R) \leq \sum_{i=1}^g {\rm deg}(\beta_i) -n \, .\] In
particular, whenever $k$ is infinite and $I$ is generically a
complete intersection, we have \[\alpha(R) \leq \sum_{i=1}^g
\delta_i -n \, .\] }
\end{itemize}
\end{remark}

\s

\begin{proposition}\label{primes}
Let $k$ be an infinite perfect field and $R$ a standard graded
$k$-algebra of dimension $d$. Let $H$ be a homogeneous $R$-ideal
such that $R_{\p}$ is regular for every minimal prime $\p$ of
dimension $d$ containing $H$. If $\alpha (R)$ is finite, there
exists a homogeneous element of degree $\alpha(R) +d$ in $0:H$ that
is not contained in any minimal prime $\p$ of dimension $d$
containing $H$.
\end{proposition}
\proof  We may assume that $\dim R/H =d$, and then by
$\p_1, \ldots, \p_s$ we denote the minimal primes of dimension
$d$ containing $H$.
Let $S$ be as in Definition~\ref{DEF} so that
$\alpha(R)=a(S)$. Write $N$ and $\P_i$ for the preimages of
$H$ and $\p_i$ in $S$. Notice that $H=NR$, that
$N_{\P_i}=0$, and that $(S/N)_{\P_i}=(R/H)_{\p_i}$ is regular for every $\P_i$.

The beginning of the proof of Proposition~\ref{syzygies}, with $g=0$,
shows that there exists a homogenous element $\beta \in 0:_S N$ of degree
$a(S)+d$ with $\beta \notin \P_i$ for each of the finitely many primes $\P_i$. Thus,
denoting the image of $\beta$ in $R$ by $\gamma$ we obtain
$\gamma \in (0:_S N)R \subset 0:_R NR = 0:_R H$ and $\gamma \notin \p_i$ for each
$\p_i$.
\QED

\s

The $\alpha$-invariant can be replaced by $a(R)$ in the above estimate
if $R$ has dimension one:

\begin{proposition}\label{dim1}
Let $k$ be an infinite field and $R$ a standard graded
Cohen-Macaulay $k$-algebra of dimension $1$. Let $H$ be a
homogeneous $R$-ideal such that $H_{\p}=0$ for every minimal prime
$\p$ of $H$. Then there exists a homogeneous element of degree $a(R)
+ 1$ in $0:H$ that is not contained in any minimal prime $\p$ of
$H$.
\end{proposition}
\proof Write $\m$ for the homogeneous maximal ideal of $R$ and set
$L=0:H$. Notice that $0: H^{\rm unm} = 0:H$ since $R \/$ is
Cohen-Macaulay. Passing to the unmixed part of $H$ we can suppose that
$R/H$ is Cohen-Macaulay of dimension $1$. Furthermore $R/L$ is
either zero or Cohen-Macaulay of dimension $1$, and $R/H + L$ has
finite length. Now the short exact sequence
\[
0 \longrightarrow R/H \cap \, L \longrightarrow R/H \oplus R/L
\longrightarrow R/H + L \longrightarrow 0
\]
induces an exact sequence of local cohomology
\[
0 \longrightarrow H^0_{\m}(R/H + L)= R/H + L \longrightarrow
H^1_{\m}(R/H \cap \, L) \, .
\]
As $a(R/H \cap \, L)\leq a(R)$ we conclude that $R/H + L$ is
concentrated in degrees $\leq a=a(R)$. Hence $\m^{a+1} \subset H +
L$. Thus any homogeneous non zerodivisor of degree $a + 1$ can
be written in the form $h + l$ where $h$, $l$ are homogeneous elements
of degree $a+1 $ in $H$, $L$ respectively. Now $l$ is an element
with the desired properties. \QED

\s

With the next two corollaries we answer the question raised at the
beginning of the section. We give a bound for the initial degrees of
annihilators of submodules generated by graded components of $\omega$
and we estimate how large the integer $i$ has to be chosen in
Corollary~\ref{omegafaithful1}(c).

\ms

\begin{theorem}\label{alpha}
Let $k$ be an infinite perfect field, let $R$ be a standard graded
reduced equidimensional $k$-algebra of dimension $d$ with
$\omega=\omega_R$, and let $t$ be an integer. If $[\omega]_{t}R$ is not a faithful $R$-module
then
\[{\rm indeg}({\rm ann}_R([\omega]_t R)) \leq \alpha(R) + d .\]
\end{theorem}
\smallskip
\proof We may assume that $\alpha(R)$ is finite. Since $R$
is generically Gorenstein, $\omega$ is isomorphic
to a suitable shift of a homogenous $R$-ideal $W$, say $\omega \cong
W(s)$. As $[W]_{t+s} R$ is not faithful it is contained in some
minimal prime $\p$ of $R$. Now
\[
{\rm ann}_R([\omega]_{t}R)={\rm ann}_R([W]_{t+s} R) \supset 0:_R \p,
\]
and the latter ideal contains a non-zero homogenous element of
degree $\alpha(R) +d$ according to Proposition~\ref{primes}. \QED

\s

Replacing Proposition~\ref{primes} by Proposition~\ref{dim1} in the proof of Theorem~\ref{alpha}
yields a better estimate for the initial degree in the one-dimensional case that
is still weaker though than the one implied by
Corollary~\ref{ann of omega dim1}.

\ms

\begin{corollary}\label{alpha2} Let $k$ be an infinite perfect field, let $R$
be a standard graded reduced Cohen-Macaulay $k$-algebra of dimension
$d \geq 1$ with homogeneous maximal ideal $\m$, and let $t$ be an integer. Write $\omega=\omega_R$
and $a=a(R)$. Then $[\omega]_{t}R$ is a faithful
$R$-module if and only if $J^{i+a+t} : \m^{a+d} \subset \m^{i}$ for
some $ i\geq \alpha(R) + d + 1$.
\end{corollary}
\proof The forward direction follows from Corollary \ref{omegafaithful1}, and the converse
from Theorems~\ref{ann of omega} and \ref{alpha}. \QED

\s

Under suitable additional assumptions the bound for $i$ in
Corollary~\ref{alpha2} can be improved, most notably by replacing
the invariant $\alpha(R)$ by the more traditional $a(R)$. This is the content of the
next three results of the section.

\ms

\begin{proposition}\label{CM}
Let $k$ be a field, let $R$ be a standard graded Cohen-Macaulay
$k$-algebra of dimension $d \geq 1$ with homogeneous maximal ideal $\m$, and
let $t$ be an integer.
Write $\omega=\omega_R$, $a=a(R)$, and
$\, -^{\vee}={\rm Hom}_R(-, \omega)$. If $\omega/(([\omega]_{\leq t}R)^{\vee \vee})$
is Cohen-Macaulay then ${\rm ann}_R([\omega]_tR)$
is a Cohen-Macaulay module with Castelnuovo-Mumford regularity $\leq d-t -1$.
In particular, whenever $t \geq -a$ one has that
$[\omega]_{t}R$ is a faithful $R$-module if and only if $J^{i+a+t} :
\m^{a+d} \subset \m^{i}$  for some $ i\geq a+d$.
\end{proposition}

\proof
In light of Corollary \ref{omegafaithful1}, the second
statement follows from the first claim and
Theorem~\ref{ann of omega}. Indeed, the bound on the Castelnuovo-Mumford regularity
shows that ${\rm ann}_R([\omega]_tR)$ is generated in degrees $\leq d-t-1 \leq a+d-1$. Hence
this annihilator vanishes by Theorem~\ref{ann of omega} if  $J^{i+a+t} :
\m^{a+d} \subset \m^{i}$  for some $i \geq a+d$.
%then Theorem~\ref{ann of omega} implies that ${\rm ann}_R([\omega]_tR) \neq R$,
%hence $t \geq -a$. Now the bound on the Castelnuovo-Mumford regularity
%shows that ${\rm ann}_R([\omega]_tR)$ is generated in degrees $\leq a+d-1$,
%and therefore vanishes by Theorem~\ref{ann of omega}.

Thus it suffices to prove the assertion about ${\rm
ann}_R([\omega]_tR)$. We may assume that this annihilator is not
zero. One has
\[{\rm ann}_R([\omega]_{t}R) = {\rm ann}_R([\omega]_{\leq t}R)={\rm ann}_R(([\omega]_{\leq t}R)^{\vee \vee}) \, . \]
The first equality is part of Theorem~\ref{ann of omega}. To see the
second equality notice that since $\omega$ is a maximal
Cohen-Macaulay $R$-module, the annihilator of every submodule is
either the unit ideal or an unmixed ideal of height zero. In
particular, two submodules of $\omega$ have the same annihilator if
they coincide locally at every minimal prime of $R$. Since
$[\omega]_{\leq t}R$ and $([\omega]_{\leq t}R)^{\vee \vee}$ are
generically equal, the second equality now follows. Thus we may
restrict our attention to the annihilator ideal ${\rm
ann}_R(([\omega]_{\leq t}R)^{\vee \vee})$. Likewise, since ${\rm
ann}_R(([\omega]_tR)^{\vee \vee}) \neq 0 ={\rm ann}_R(\omega)$ it follows that
$([\omega]_{\leq t}R)^{\vee \vee}$ and $\omega$ cannot coincide
locally at every minimal prime of $R$. Thus $\omega/(([\omega]_{\leq
t}R)^{\vee \vee})$ has dimension $d$, hence is a maximal
Cohen-Macaulay $R$-module.

Write $R=S/I$ with $S=k[X_1, \ldots,X_n]$ a polynomial ring and $I$ a homogeneous $S$-ideal of height
$g$. Choose a regular sequence $\ul{\beta}=\beta_1, \ldots, \beta_g$ of
forms of degree $\delta \gg 0$ contained in $I$ and set $L= (\ul{\beta}) :_S I$.
Notice that $S/L$ is a $d$-dimensional Cohen-Macaulay ring. Moreover one has
$(L/(\ul{\beta}))(g \delta - n) \cong \omega$. Thus there exists a homogeneous ideal
$H$ of $S$ so that $(\ul{\beta}) \subset H \subset L$ and
$(H/(\ul{\beta}))(g \delta -n) \cong ([\omega]_{\leq t}R)^{\vee \vee}$. Clearly
\[{\rm ann}_R(([\omega]_{\leq t}R)^{\vee \vee}) = ((\ul{\beta}):_S H)/I \, . \]

Since $L/H$ is isomorphic to a shift of the module
$\omega/(([\omega]_{\leq t}R)^{\vee \vee})$ it follows that $L/H$
is a maximal Cohen-Macaulay $R$-module. Therefore $S/H$ is a
$d$-dimensional Cohen-Macaulay ring, and hence so is
$S/((\ul{\beta}):_SH)$. Thus $((\ul{\beta}):_S H)/I$ is a maximal
Cohen-Macaulay $R$-module, proving the assertion in the proposition
about the Cohen-Macaulayness of the annihilator ideal.

Applying ${\rm Ext}^g_S(-,S(-n))$ to the exact sequence
\[ 0 \longrightarrow ((\ul{\beta}):_SH)/I \longrightarrow S/I
\longrightarrow S/((\ul{\beta}):_S H) \longrightarrow 0
\]
of maximal Cohen-Macaulay $R$-modules we obtain this exact sequence,
\[ 0 \longrightarrow \omega_{S/((\ul{\beta}): H)} \longrightarrow \omega_{S/I}=\omega
\longrightarrow {\rm{Ext}}^g_S(((\ul{\beta}):_SH)/I, S(-n)) \longrightarrow 0 \, .
\]
The canonical module on the left satisfies
\[\omega_{S/((\ul{\beta}): H)} \cong (H/(\ul{\beta}))(g\delta -n)
\cong ([\omega]_{\leq t} R)^{\vee \vee} \supset [\omega]_{\leq t}R \, , \]
where the first isomorphism holds because $H$ is unmixed.
Therefore the above exact sequence shows that
${\rm{Ext}}^g_S(((\ul{\beta}):_SH)/I, S(-n))$ is concentrated in
degrees $\geq t+1$, and hence by local duality the local cohomology module
$H^d_{\m}(((\ul{\beta}):_SH)/I)$ is concentrated in degrees $\leq -t-1$. As
$((\ul{\beta}):_SH)/I$ is a $d$-dimensional Cohen-Macaulay module we deduce that
it has regularity $\leq d-t-1$. \QED

\s

We follow suit with an estimate for rings of type $2$ that
does not require the Cohen-Macaulayness of $\omega/(([\omega]_{\leq t}R)^{\vee \vee})$.

\ms

\begin{proposition}\label{type2}
Let $k$ be a field and $R$ a standard graded geometrically reduced
Cohen-Macaulay $k$-algebra of type $2$ and dimension $d \geq 1$ with
homogeneous maximal ideal $\m$. Write $R=S/I$ with $S=k[X_1,
\ldots, X_n]$ a polynomial ring and $I$ a homogeneous $S$-ideal of
height $g$. Consider the last map in a minimal homogeneous free
$S$-resolution of $R$
\[
0 \longrightarrow S(-l_1) \oplus S(-l_2)
\stackrel{\varphi}\longrightarrow \oplus_{i=1}^{m} S(-k_i) \, ,
\]
where $l_1 \leq l_2$ and $k_1 \leq \ldots \leq k_m$. Write
$\omega=\omega_R$ and $a=a(R)$. If $[\omega]_{-a}R$ is not a faithful
$R$-module then
\[{\rm indeg}({\rm ann}_R([\omega]_{-a}R)) \leq gl_1 +l_2 -
\sum_{i=1}^{g+1} k_i -g =a+d +gl_1 - \sum_{i=1}^{g+1} k_i\, .\]
Equivalently, $[\omega]_{-a}R$ is a faithful $R$-module if and only
if $J^{i} : \m^{a+d} = \m^{i}$  for some  $ i\geq a+d +gl_1 -
\sum_{i=1}^{g+1} k_i+1$.
\end{proposition}
\proof  By Theorem~\ref{ann of omega} and Corollary \ref{omegafaithful2} it suffices to prove the first
statement. We may assume that the field $k$ is infinite and perfect. Since $R$
is Cohen-Macaulay we obtain the following presentation,
\[
 \oplus_{i=1}^{m} S(k_i) \stackrel{\varphi^*}\longrightarrow S(l_1)
\oplus S(l_2)\longrightarrow \omega(n) \longrightarrow 0 \,.
\]
Using the standard bases of the free modules this presentation gives
rise to homogeneous generators $w_1, w_2$ of $\omega(n)$ and a
matrix
\[
\left(
               \begin{array}{ccc}
                 f_1 & \cdots  & f_m \\
                 h_1 & \cdots  & h_m \\
               \end{array}
             \right)
\]
representing $\varphi^*$. Notice that ${\rm deg}(w_1)=-l_1 \geq
{\rm deg}(w_2)=-l_2$, ${\rm deg}(f_i)=l_1-k_i$, and ${\rm deg}
(h_i)=l_2-k_i$. One has ${\rm deg} (w_1) > {\rm deg}(w_2)$ since
otherwise $[\omega]_{-a}R= \omega$ would be faithful. Therefore \[{\rm
ann}_R ([\omega]_{-a} R) = {\rm ann}_R \, w_2 \,.\] Further observe that
\[w_2R :_S w_1= (f_1, \ldots, f_m)\]
and
\begin{equation}\label{omegasyz}
{\rm ann}_S \, w_2 = \displaystyle{\{} \sum \lambda_ih_i \ | \  \sum
\lambda_i f_i =0 \}\, ;
\end{equation}
the last equality obtains because an element $\varepsilon$ of $S$
belongs to ${\rm ann}_S \, w_2 $ if and only if the vector $\left(
                                                              \begin{array}{c}
                                                                0 \\
                                                                \varepsilon \\
                                                              \end{array}
                                                            \right)$
is in the column space of the above matrix.

%One has \[I \subset w_2R :_S w_1=(f_1, \ldots, f_m) \, .\] \noindent
As ${\rm ann}_R \, w_2 \not=0$  and $R$ is
unmixed, there exists a minimal prime $\p$ of $R$ such that $({\rm
ann}_R \, w_2)_{\p} \not=0$. In particular, $w_2R_{\p} \not= w_1
R_{\p} + w_2 R_{\p}$ because the latter module is faithful. Thus the ideal
$(f_1, \ldots, f_m) =w_2R :_S w_1$ is contained in $\P$, the
preimage of $\p$ in $S$. On the other hand, this ideal contains $I$. It follows
that $(f_1, \ldots, f_m)$ has height $g$. Furthermore,
since $I$ is radical the localization
$(f_1, \ldots , f_m) S_{\Q}$ is a complete intersection prime ideal
for every prime $\Q$ of height $g$ in $S$ containing $(f_1, \ldots, f_m)$.
Finally, $g=\height \,I=\height \,I_2(\varphi)
\leq m -1$ by the Eagon-Northcott bound on the height of
determinantal ideals. Now Proposition~\ref{syzygies} shows that
locally at $\P$, the syzygy module of $(f_1, \ldots, f_m)$ is
generated by its elements of degrees $\leq \sum_{i=1}^{g+1} (l_1-
k_i) + n - g + a(S)$ and by the Koszul relations. Therefore \ref{omegasyz}
gives
\[
({\rm ann}_S \, w_2)_{\P}=([{\rm ann}_S \, w_2]_{\leq
\sum_{i=1}^{g+1} (l_1- k_i) + n - g + a(S) + l_2 -l_1} +
I_2(\varphi^*))_{\P} \, .
\]
As $I_2(\varphi^*) \subset {\rm ann}_S \, \omega =I$ we conclude
that
\begin{eqnarray*}
({\rm ann}_R \, w_2)_{\p}&=&([{\rm ann}_R \, w_2]_{\leq
\sum_{i=1}^{g+1} (l_1- k_i) + n - g + a(S) + l_2 -l_1})_{\p}\\
&=&([{\rm
ann}_R \, w_2]_{\leq g  l_1 +l_2 - \sum_{i=1}^{g+1} k_i -g})_{\p} \, .
\end{eqnarray*}
The assertion now follows since $({\rm ann}_R \,
w_2)_{\p}\not=0 \, $.\QED

\ms

\begin{corollary}\label{codim2}
Let $k$ be a field and let $R$ be a standard graded equidimensional
geometrically reduced $k$-algebra of dimension $d \geq 1$ with homogeneous
maximal ideal $\m$. Assume that $R$ is an almost complete intersection
of embedding codimension $2$.  Write $\omega=\omega_R$ and $a=a(R)$. If
$[\omega]_{-a}R$ is not a faithful
$R$-module then
\[{\rm indeg}({\rm ann}_R([\omega]_{-a}R)) \leq a + d \, .\]
Equivalently, $[\omega]_{-a}R$ is a faithful $R$-module if and only
if $J^{i} : \m^{a+d} = \m^{i}$  for some  $ i\geq a+d+1$.
\end{corollary}
\proof Notice that $R$ is Cohen-Macaulay by the Syzygy Theorem and has
type $2$, see \cite[2.1]{EG}. Hence we may apply
Proposition~\ref{type2} with $g=2$. It suffices to show that $2l_1 -
\sum_{i=1}^{3} k_i \leq 0$. This holds, because $k_3=l_1+l_2
-k_1-k_2$ by the Hilbert-Burch Theorem and $l_1 \leq l_2 \, $. \QED

\s

We finish this section with a different estimate for initial degrees of
annihilators -- an estimate from below. In the proof we use the notation
$H^{\rm unm}$ for the {\it unmixed part} of an ideal $H$, which is the intersection
of the primary components of maximal dimension.

\ms

\begin{proposition}
Let $k$ be a field, let $R$ be a standard graded
Cohen-Macaulay $k$-algebra of dimension $d$, and let $H $ be
a homogeneous $R$-ideal. Write $c=c(R)$ as in $\ref{standard}$.
One has  ${\rm indeg} (0:H) \geq c +d +1 -
e(R/H) $.
\end{proposition}
\proof We may assume that $k$ is infinite and that $0:H \neq 0$.
Set $e= e(R/H)$. We prove the claim by induction on $d$.
First let $d=0$. In this case $e = \lambda (R/H)$. Therefore
$\m^e \subset H$, which gives $0: H \subset 0 :\m^e$. But
${\rm indeg} (0:\m^e) \geq c +1 - e$ because
$c$ is the initial degree of the
socle of $R$.

%Choose a homogeneous element $f \not= 0$ in $0:H$. As
%$0: f \supset H$ we obtain this inequality of lengths,
%\[
%\lambda(fR) =\lambda (R/(0:f)) \leq \lambda (R/H)= e \, .
%\]
%Therefore $\m^e f=0$ and we deduce that $e-1 + {\rm deg} (f) \geq c \, $.
Next let $d \geq 1$. Notice that $0: H^{\rm unm} = 0:H$ since $R \/$ is
Cohen-Macaulay and that $e(R/H^{\rm unm})=e(R/H)$ by the associativity
formula for the multiplicity. Hence we may replace $H$ by $H^{\rm
unm}$ to assume that $H$ is unmixed. Since $0:H
\not= 0$ we have $\height \, H =0$. Furthermore, $0:H$ is unmixed
with $\height  (0:H) =0$ or else $0:H = R$. Now let $x \in R$ be a
linear form that is regular on $R$, and write $^{^{\tratto}}$ for images in
the $d-1$ dimensional Cohen-Macaulay ring $\ol{R}=R/(x)$. Notice
that $x$ is regular modulo $H$ and is regular modulo $0:H$ unless
$0:H=R$. It follows that  $e(\ol{R}/\ol{H})=e(R/H)$ and ${\rm
indeg}(\ol{0:_R H}) ={\rm indeg}(0 :_R H)$. Furthermore,
$c(\ol{R})=c(R) +1$ as $\omega_{\ol{R}} \cong \omega_R (1)$. Thus we
conclude
\begin{eqnarray*}
  {\rm indeg} (0 :_R H) &=& {\rm indeg}( \ol{0:_R H})  \\
   &\geq& {\rm indeg}(\ol{0}:_{\ol{R}} \ol{H})  \hspace{3.765cm} \hbox{\rm since \ $ \ol{0:_R H} \subset \ol{0}:_{\ol{R}} \ol{H} $}\\
   &\geq& c(\ol{R}) + {\rm dim} \, \ol{R} +1 -  e(\ol{R}/\ol{H})  \hspace{1.6cm} \hbox{\rm by induction hypothesis}\\
   &=& (c +1) + (d-1) +1 -e\\
   &=& c+d+1 -e \,.
\end{eqnarray*} \QED

\s

If $R$ is zero-dimensional in the above proposition, then the lower bound $c+1-e(R/H)$ can be improved to $c-a(R/H)$.
This shaper estimate immediately gives the inclusion $J^i:\m^j \subset \m^{i-j+c+d}$ in the
setting of Corollary \ref{colonmax1}.

\s

\section{The core of standard graded algebras}

\sss

In this section we apply the previous results to the core of powers of
the maximal ideal.

\begin{Assumptions}\label{asscore}
{\rm In addition to the assumptions of \ref{standard} suppose that
if ${\rm char} \,  k > 0$  then $k$ is infinite and $R$ is
geometrically reduced.}
\end{Assumptions}

\ms

The next result is essentially present, in a more general form, in
\cite[4.2]{FPU}.

\begin{theorem}\label{formula}  With assumptions as in $\ref{asscore}$ one has for every $n \geq 1$,
\[
\core{\m^n} = J^{nd+a+1} : \m^{a+d}.
\]
\end{theorem}
\proof From \cite[2.3 and 2.5]{PUV} we know that \[\core{\m^n}
=(J^{[n]})^{j+1}:(\m^n)^j \ \ \  \hbox{\rm for} \ \  j \gg 0 \, .\]
Furthermore
\begin{eqnarray*}
\nonumber %to remove numbering (before each equation)
 (J^{[n]})^{j+1}:(\m^n)^j&=& (J^{[n]})^{[j+1]}:(\m^n)^{jd}   \hspace{1.6cm} \hbox{\rm by (\ref{PUV 2.2})} \\
  &=& J^{[nj+n]}:\m^{njd} \\
  &=& J^{nj+n}: \m^{nj+n-nd+d-1}  \hspace{0.951cm} \hbox{\rm by (\ref{PUV 2.2})} \\
  &=& J^{nd+a+1} : \m^{a+d}  \hspace{2.055cm} \hbox{\rm by
  Remark~\ref{colon1}.}
\end{eqnarray*}\QED

\s

The above theorem relates the core of powers of the maximal ideal to
the colon ideals studied in the previous sections. We leave it to
the reader to express most of the earlier results in terms of cores.
Here we only collect the main applications:

\ms

\begin{corollary}\label{coreandK}
With assumptions as in $\ref{asscore}$ one has for every $n \geq 1$ $:$
\begin{itemize}
\item [(a)] $\core{\m^n}=\m^{nd+a+1}+N$ for some ideal $N$ of height zero

\item[(b)] $\m^{nd+a+1} \subset \core{\m^n} \subset \m^{nd+b+1}$.
\end{itemize}
\end{corollary}
\proof The result follows immediately from Theorem~\ref{formula} and
Corollaries~\ref{K}
 and \ref{colonmax1}. \QED

\s

Item (b) of the next theorem was asserted in \cite[4.1]{HySm2},
assuming only that $R$ is reduced. However, as the theorem shows, this
statement is equivalent to the faithfulness of the module
$[\omega]_{-a}R$. Furthermore, $R$ being {\it geometrically} reduced
is essential according to \cite[5.1]{FPU}.

\ms

\begin{theorem}\label{core-ann} With assumptions as in $\ref{asscore}$ the following
are equivalent $:$
\begin{itemize}
  \item [(a)] $[\omega]_{-a}R$ is faithful
  \item [(b)] $\core{\m^n}=\m^{nd+a+1}$ for every $n\geq 1$
  \item [(c)] $\core{\m^n}$ is generated in one degree for some $n
  \gg 0 \, $.
\end{itemize}
\end{theorem}
\proof Theorem~\ref{formula} and Corollary~\ref{omegafaithful2} show
that (a) implies (b). If (c) holds then according to
Theorem~\ref{formula} and Corollary~\ref{power} one has $ J^{nd+a+1}
: \m^{a+d}=\m^{nd+a+1}$. Notice that $nd+a+1 \gg 0$. Thus again
Corollary~\ref{omegafaithful2} gives that (a) obtains. \QED

\ms

\begin{corollary}\label{indeg=a+d}
In addition to the assumptions of $\, \ref{asscore}$ suppose that one of
the following conditions holds $:$
\begin{itemize}
  \item $R$ is reduced and  $ [\omega]_{\geq d + 1}=\core \m  \, \omega$
  \item $d=1$
  \item $R$ is a reduced almost complete intersection of embedding
  codimension $2$.
\end{itemize}
Then the following are equivalent $:$
\begin{itemize}
  \item [(a)] $[\omega]_{-a}R$ is faithful
  \item [(b)] $\core{\m}=\m^{a+d+1}$
  \item [(c)] $\core{\m}$ is generated in one degree.
\end{itemize}
\end{corollary}
\proof We apply Theorem~\ref{formula} and Corollary~\ref{power}. In the first case we also use
Corollary~\ref{omegafaithful2} via Theorem~\ref{omega}(a), in the second case
Corollary~ \ref{omegafaithfuldim1}, and in the third case Corollary~\ref{codim2}. \QED

\s

\section{The core of points}

\sss

\begin{Assumptions and Discussion}\label{core of points}
{\rm Let $k$ be an infinite field and let $X=\{P_{1}, \ldots,
P_{s}\}$ be a set of $s$ reduced points in $\mathbb{P}_{k}^{n}$.
Write $S=k[x_0, \ldots, x_n]$ for the polynomial ring and $R=S/I_X$
for the homogeneous coordinate ring of $X \subset
\mathbb{P}_{k}^{n}$. Let $\m$ denote the homogeneous maximal ideal
of $R$, $K$ its total ring of quotients, $B$ the integral closure of
$R$ in $K$, and $\C= R :_K B$ the conductor. Furthermore, we write
$\omega=\omega_R$, $a=a(R)$, $b=b(R)$ and we define
$\core{X}=\core{\m}$. Finally, let $y \in R$ be a linear form that
is $R$-regular and set $J=yR$. Notice that $a \leq s-2$ and that $R$
is geometrically reduced.

Homogeneous polynomials $f_{1}, \ldots, f_{s}$ in $S$ are called
{\it separators} of $X$ if $f_{i}(P_{j})=\delta_{ij}$ for every $i$,
$j$. They are called {\it minimal separators} if in addition each
$f_i$ has smallest possible degree.}

\end{Assumptions and Discussion}
\ms
With the next lemma we recall a known fact describing the
conductor in terms of minimal separators (see for instance
\cite[3.13]{GKR}). We include a proof for the convenience of the
reader.
\ms

\begin{lemma}\label{conductor}
With assumptions as in $\ref{core of points}$ let $h_{1}, \ldots,
h_{s}$ be a collection of separators of $\/X$ and $f_{1}, \ldots,
f_{s}$ a collection of minimal separators. One has $h_i R \subset
f_i R \/ $ and $f_{1}, \ldots, f_s$ minimally generate $\C$ as an
$R$-ideal; in particular,   $\m^{a+1} \subset (f_1, \ldots, f_s)R =
\C$ and ${\rm deg} (f_i) \leq a+1$.
\end{lemma}
\proof
Fix $\lambda_{ij} \in k$ with $P_{i}=(\lambda_{i_0} : \ldots : \lambda_{i_n})$.
We use these identifications of non-negatively graded rings,
\[\varphi: R \hookrightarrow B \cong k[t] \times \ldots \times  k[t] = \oplus_{i=1}^{s} k[t]e_i \, ,\]
where $k[t]$ is a standard graded polynomial ring, $e_i$ are
standard basis elements of degree zero, and $\varphi$ maps the image
in $R$ of a polynomial $h\in [S]_l$ to the tuple
$\sum h(\lambda_{i_0}, \ldots ,\lambda_{i_n}) t^{l} e_i$. Thus the
elements $t^{l_1}e_1, \ldots, t^{l_s}e_s$ belong to $R$ if and only
if there exist separators $h_1, \ldots , h_s$ of degrees $l_1, \ldots,
l_s$, in which case $t^{l_i}e_iR=h_iR$. From this we see that $h_i R
\subset f_i R \/ $. It also follows that the minimal separators
$f_{1}, \ldots, f_{s}$ minimally generate the largest $R$-ideal of
the form $\sum t^{l_i} e_i R= \sum t^{l_i} k[t] e_i$, equivalently,
the largest homogeneous $B$-ideal contained in $R$. However, this
ideal is the conductor $\C$.

Finally, the long exact sequence of local cohomology shows that
$B/R$ is concentrated in degrees $\leq a$, hence $\m^{a+1} \subset
\C$. \QED

\s

The next proposition gives a geometric interpretation of the core of
points in terms of separators:

\ms

\begin{proposition}\label{separators}With assumptions as in $\ref{core of points}$
let $f_{1}, \ldots, f_{s}$ be minimal separators of $X$.
One has
\begin{itemize}
  \item [(a)] $\core{X}=y\C=\m\C =y R (f_{1}, \ldots, f_{s})=\m(f_{1}, \ldots, f_{s})$
  \item [(b)]  $ [\omega]_{\geq 2}=\core X \, \omega$.
\end{itemize}
\end{proposition}
\proof According to Theorem~\ref{formula} one has $\core X= J^{a+2}
: \m^{a+1}$. Now the first and second equality in (a) follow from
Remark~\ref{coreandS}, and (b) is a consequence of
Proposition~\ref{=indim1}. Finally, Lemma~\ref{conductor} implies
$\C=(f_{1}, \ldots, f_{s})R$. \QED

\s

\begin{corollary}\label{Y and Z}
In addition to the assumptions of $\, \ref{core of
points}$ suppose that $X=Y \cup Z$, where $Y$ is contained in a
hypersurface $f=0$ and $Z$ is a collection of $e$ reduced points whose
homogeneous coordinate ring has $a$-invariant $a'$. One has \[
\m^{a+2} + f \m^e \subset \m^{a+2} + f \m^{a' + 2} \subset \core X
\subset \m^{b+2} \, .\]
%in particular $b \leq a' + {\rm deg}(f)$.
\end{corollary}
\smallskip
\proof  Let $h_1, \ldots, h_e$ be minimal separators of $Z$, and
write $H$ for the defining ideal of $Z$ in $X$. From
Lemma~\ref{conductor} we know that $\m^{a'+1} \subset (h_1, \ldots,
h_e) R + H$. Since $fH=0$ in $R$, multiplying this equation by $f$ we
obtain $f\m^{a'+1} \subset (fh_1, \ldots , fh_e)R$. However, those elements of $fh_1,
\ldots , fh_e$ that are not contained in $I_{X}$ form part of a collection of separators of $X$.
Hence $(fh_1, \ldots , fh_e)R \subset \C$ by the same Lemma~\ref{conductor}, and therefore
$f \m^{a'+2} \subset \core X$ according to
Proposition~\ref{separators}(a). The remaining assertions follow from
Discussion~\ref{core of points} and Corollary~\ref{coreandK}(b). \QED

\s

The previous result suggests that the shape of the core is related to
uniformity properties of the set of points. We recall one such
condition: The scheme $X \subset \mathbb{P}_{k}^{n}$ is said to have
the {\it Cayley-Bacharach} property if each subscheme of the form $X
\backslash \{P_i\} \subset \mathbb{P}_{k}^{n}$ has the same Hilbert
function.

\ms

\begin{corollary}\label{CB} With assumptions as in $\ref{core of points}$
the scheme $X$ has the Cayley-Bacharach property if and only if
$\core X =\m^{a+2}.$
\end{corollary}
\proof It is easy to see that $X$ has the Cayley-Bacharach property
if and only if the minimal separators all have the same degree.
According to Proposition~\ref{separators}(a) this means that
$\core X$ is generated in one
degree, which in turn is equivalent to $\core X =\m^{a+2}$ as shown
in Corollary~\ref{indeg=a+d}.

Alternatively, in \cite[3.5]{GKR} the Cayley-Bacharach property has
been characterized in terms of the faithfulness of $[\omega]_{-a}R$.
Again according to Corollary~\ref{indeg=a+d} the latter condition
holds if and only if $\core X =\m^{a+2}$. \QED

\s

The next example illustrates the previous two corollaries.

\ms

\begin{example}
{\rm In addition to the assumptions of \ref{core of points} suppose
${\rm char} \, k \not=2$ and take $X$ to be the $4$ points
$(0:-1:1), (0:0:1), (0:1:1), (1:0:1)$ in $\mathbb{P}^2_k$. These
points and their separators are depicted in this figure:

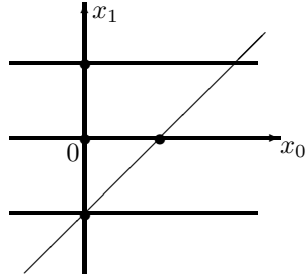
\setlength{\unitlength}{1cm}
\begin{picture}(6,4)

%\showgrid
%\set{xunit=1cm, yunit=1cm}

\put(2,2){\vector (1,0){3.6}}

\put(3,.2){\vector (0,1){3.6}}

\put(5.6,1.8){$x_0$}

\put(3.1,3.6){$x_1$}

\put(2.92,1.89){$\bullet$} \put(2.76,1.69){$0$}

\put(2.92,2.89){$\bullet$}

\put(2.92,0.89){$\bullet$}

\put(3.92,1.89){$\bullet$}

\linethickness{.04cm}

\put(2,2){\line (1,0){3.6}}

\put(3,.2){\line(0,1){3.6}}

\put(2,3){\line (1,0){3.3}}

\put(2,1){\line (1,0){3.3}}

\put(2.2,0.2){\line (1,1){3.2}}

\end{picture}

\noi Notice that
\[
R=k[x_0,x_1,x_2]/(x_0 x_1,x_0(x_0-x_2),x_1(x_1-x_2)(x_1+x_2)),
\]
and one easily sees that $a=1$. We choose $y$ to be the image of
$x_2$, and as minimal separators of $X$ we take $x_1(x_1-x_2),
(x_0-x_1-x_2)(x_1-x_2), x_1(x_1+x_2), x_0$. From the geometric
interpretation of the core in terms of separators,
Proposition~\ref{separators}(a), one immediately sees that $\core X=
\m^3 + x_0^2R \supsetneq \m^3=\m^{a+2}$. The inclusion $x_0^2R
\subset \core X$ would have also been predicted by Corollary~\ref{Y
and Z} with $f=x_0$, and the strict containment $\m^{3} \subsetneq
\core X$ reflects the obvious fact that $X$ does not have the
Cayley-Bacharach property. }
\end{example}

\ms

\begin{corollary} Let $k$ be a field of characteristic zero, let
$Y\subset \mathbb{P}^{n+1}_k$ be a reduced and irreducible
arithmetically Cohen-Macaulay curve, and write $\n$ for the
homogeneous maximal ideal of the homogeneous coordinate ring of $Y$.
Consider a general hyperplane section $X\subset \mathbb{P}^{n}_k$ of
$Y$ and use the notation of $\, \ref{core of points}$. One has $\core Y=
\core \n= \n^{a + 2}$ and $\core X= \m^{a +2}$.
\end{corollary}
\proof According to \cite[3.4]{H} the set of points $X$ has the
 Cayley-Bacharach property. Now the two equalities follow
 from Corollaries~\ref{coreandK} and  \ref{CB}.\QED

\s
\section{Local estimates on cores}
\sss

We finish this paper with a generalization of Corollary~\ref{Y and
Z} to the context of zero-dimensional ideals in local rings.
\ms

\begin{proposition}\label{local}
Let $(R, \m)$ be a local Cohen-Macaulay ring with infinite residue
field, $I$ an $\m$-primary $R$-ideal, $L$ and $H$ two $R$-ideals
such that $LH=0$ in $R$. Write $e=e(I; R/H)$ for the multiplicity of
the ring $R/H$ with respect to the ideal $I$. Then $L I^{e} \subset
\core I$.
\end{proposition}
\proof We prove the claim by induction on $d ={\rm dim}\, R$. If
$d=0$ then $e$ is the length of $R/H$, and we easily obtain $I^e
\subset H$. Therefore $L I^e \subset L H =0$. Now consider the case
$d \geq 1$. For $J$ an arbitrary reduction of $I$ we need to show that $L
I^{e} \subset J$. Let $x$ be a general element of $J$ and
write $^{^{\tratto}}$ for images in the $d-1$
dimensional Cohen-Macaulay ring $\ol{R}=R/(x)$. Notice that $L
H^{\rm unm}=0$ since $R$ is Cohen-Macaulay and that $e(I; R/H^{\rm
unm})=e$ by the associativity formula for the Hilbert-Samuel
multiplicity. Thus we may replace $H$ by $H^{\rm unm}$ to assume
that $H$ is unmixed. We may further suppose that $L \not= 0$.
Therefore ${\rm dim}\, R/H= {\rm dim}\, R=d$ since $R$ is
Cohen-Macaulay. We conclude that the general element $x$ of $J$ is
regular on $R/H$. Thus $e(\ol{J}; \ol{R}/\ol{H})=e(J; R/H)$. As
$\ol{J}$ and $J$ are reductions of $\ol{I}$ and $I$, respectively,
we have $e(\ol{J}; \ol{R}/\ol{H})=e(\ol{I}; \ol{R}/\ol{H})$ and
$e(J; R/H)=e(I; R/H)$. It follows that $e(\ol{I}; \ol{R}/\ol{H})=e(I; R/H)=e$.
Now our induction hypothesis gives $\ol{L I^e} \subset \core{\ol{I}}
\subset \ol{J}$. Hence indeed $L I^{e}
\subset J$. \QED

\s

We obtain the estimate $f\m^e \subset \core X$ of
Corollary~\ref{Y and Z} from the above proposition if we take $L=fR$
and $H= I_{Z}R$.

\s
\s

\s
\s

\end{document}